\documentclass[12pt,draft]{article}

\usepackage[T1]{fontenc}
\usepackage{amsfonts}

\newtheorem{thm}{theorem}[section]
\newtheorem{theorem}[thm]{Theorem}
\newtheorem{proposition}[thm]{Proposition}
\newtheorem{lemma}[thm]{Lemma}
\newtheorem{corollary}[thm]{Corollary}
\newtheorem{remark}[thm]{Remark}

\newtheorem{definition}[thm]{Definition}

\begin{document}

\title{On the Graded Identities for Elementary Gradings in Matrix Algebras over Infinite Fields}

\author{Diogo Diniz Pereira da Silva e Silva \footnote{Supported by CNPq}}
\date{{\small \textit{Unidade Acadêmica de Matemática e Estatística\\
Universidade Federal de Campina Grande\\
Cx. P. 10.044, 58429-970, Campina Grande, PB, Brazil}}\\
E-mail: diogo@dme.ufcg.edu.br}
\maketitle

\begin{abstract}
We find a basis for the $G$-graded identities of the $n\times n$ matrix algebra $M_n(K)$ over an infinite field $K$ of characteristic $p>0$ with an elementary grading such that the neutral component corresponds to the diagonal of $M_n(K)$.
\end{abstract}

\section{Introduction}

The polynomial identities of the matrix algebra $M_n(K)$ are important in the theory of PI-algebras, for example, the T-ideal of its polynomial identities $T(M_n(K))$ arises in Kemer's structural theory of T-ideals as one of the  T-prime T-ideals. However, over infinite fields, finite bases for $T(M_n(K))$ were determined only when $n=2$ and $char$ $K \neq 2$. In (\cite{Razmyslov}) Razmyslov  determined a basis of the identities of $M_2(K)$ with $9$ elements in the case $char$ $K =0$ and in (\cite{Drensky}) Drensky improved this result by finding a minimal basis with two identities. A basis for the identities of $M_2(K)$ over infinite fields of characteristic $p>2$ was determined by Koshlukov in (\cite{Koshlukov}). He has proved that the same basis found in (\cite{Drensky}) is a basis for the identities of $M_2(K)$ over an infinite field $K$ of characteristic $>3$, and if $char$ $K =3$ one more identity is necessary. The problem of determining a finite basis for $T(M_n(K))$ when $n\geq 3$ is still open.

In \cite{Kemer} Kemer develops a theory of T-ideals analogous to the theory of ideals in commutative polynomial algebras and the concept of $\mathbb{Z}_2$-graded identities was a key component in this theory showing the importance of $\mathbb{Z}_2$-graded identities. In \cite{Vincenzo} Di Vincenzo found a basis for the $\mathbb{Z}_2$-graded identities of $M_2(K)$, over a field of characteristic $0$.

Soon afterwards, the study of $G$-graded identities of algebras graded by an arbitrary group $G$ became a problem of independent interest. Vasilovsky \cite{Vasilovsky1, Vasilovsky2} extended the results of \cite{Vincenzo} and determined basis for the graded identities of $M_n(K)$ graded by the groups $\mathbb{Z}$ and $\mathbb{Z}_n$, for any $n$, in the case $K$ is a field of characteristic $0$. Azevedo \cite{Azevedo1, Azevedo2} proved that the results of Vasilovsky also hold if $K$ is an infinite field of characteristic $p>0$. 

The complete description of all possible $G$-gradings of $M_n(K)$ by a finite group $G$, over an algebraically closed field $K$ of characteristic $0$, is given in \cite{Bahturin}. The main result states that any such grading is the tensor product of two types of gradings: the \textit{fine gradings} where $\dim (M_n(K))_g \leq 1$ for all $g \in G$ and the \textit{elementary gradings} which are induced by a grading in the vector space $K^n$. 
For the elementary gradings the $G$-graded identities of $M_n(K)$, where $G$ is an arbitrary group, were studied in \cite{drenskyebahturin} in a very general setting. Provided that the neutral component coincides with the diagonal of $M_n(K)$ a basis of the $G$-graded identities was determined when $char K =0$.

In this article we combine the methods of \cite{drenskyebahturin} and \cite{Azevedo1, Azevedo2} to prove that the same results of \cite{drenskyebahturin} for the graded identities of $M_n(K)$ with elementary gradings hold for infinite fields of characteristic $p>0$.

\section{Preliminaries}\label{preliminaries}
In this article $K$ denotes an infinite field and all vector spaces and algebras will be considered over $K$. Let $G$ be an arbitrary group, a \textit{$G$-grading} of an algebra $A$ is a vector space decomposition $A=\bigoplus_{g \in G} A_g$ such that for any $g,h \in G$ the inclusion $A_gA_h\subset A_{gh}$ holds. The elements of $A_g$ are said to be \textit{homogeneous of degree $g$}. We denote by $\epsilon$ the identity element of $G$, the component $A_{\epsilon}$ is called the \textit{neutral component}.

We denote by $I_k$ the set $\{1,2,\dots, k\} \subset \mathbb{N}$, given $i,j \in I_n$ we denote by $E_{ij} \in M_n(K)$ the matrix unit in which the only non-zero entry is $1$ in the $i$-th row and $j$-th column. Given an 
$n$-tuple $(g_1,\dots, g_n) \in G^{n}$ a $G$-grading in $M_n(K)$ is determined by imposing that $E_{ij}$ is homogeneous of degree $g_i^{-1}g_j$. These gradings are called \textit{elementary gradings}. 

\begin{proposition}[\cite{IoneMontes}]
Let $G$ be a group, the $G$-grading of $M_n(K)$ is elementary if and
only if all matrix units $E_{ij}$ are homogeneous.
\end{proposition}

Let $\{X_g| g\in G\}$ be a family of disjoint countable sets indexed by $G$ and let $X=\cup_{g \in G}X_g$. We denote by $K\langle X \rangle$ the free associative algebra freely generated by $X$. Given a monomial $m=x_{i_1}\dots x_{i_k}$ we denote by $h(m)$ the $k$-tuple $(h_1,\dots, h_k)$ where for each $l \in I_k$ the $l$-th coordinate is the element $h_l\in G$ such that $x_{i_l}\in X_{h_l}$. We denote by $K\langle X \rangle_{g}$ the subspace of $K\langle X \rangle$ generated by the monomials $x_{j_1}\dots x_{j_k}$ such that $h_1\cdot\dots \cdot h_k=g$, where $(h_1,\dots,h_k)=h(m)$. The decomposition $K\langle X \rangle = \bigoplus_{g \in G} K\langle X \rangle_{g}$ is a $G$-grading and with this grading $K\langle X \rangle$ is the free $G$-graded associative algebra on $X$. In this grading any monomial $m=x_{i_1}\dots x_{i_k}$ is homogeneous and its degree with respect to the $G$-grading will be denoted by $\alpha(m)$.

A polynomial $f(x_1,\dots,x_k) \in K\langle X \rangle$ is a \textit{graded polynomial identity} for the $G$-graded algebra $A$ if $f(a_{1},\dots, a_{k})=0$ whenever $a_l \in A_{h(x_l)}$ for every $l \in I_k$. We denote by $T_G(A)$ the set of all graded identities of the $G$-graded algebra $A$, this set is an ideal of $K\langle X \rangle$ which is invariant under all graded endomorphisms of $K\langle X \rangle$. It is easy to show that the intersection of a family of $T_G$-ideals of $K\langle X \rangle$ is also a $T_G$-ideal, hence given $S\subset K\langle X \rangle$ we may define the $T_G$-ideal generated by $S$, denoted by $\langle S \rangle^{T_G}$, as the intersection of all $T_G$-ideals that contain $S$. We say that $S$ is a basis of the graded identities of $A$ if $T_G(A)=\langle S \rangle^{T_G}$.

In this article, with the exception of Proposition \ref{p1} and Remark \ref{r1}, we fix an $n$-tuple $\textbf{g}=(g_1,\dots, g_n)\in G^{n}$ of pairwise different elements and let $M_n(K)=\bigoplus_{g \in G}(M_n(K))_g$ be the elementary grading induced by $\textbf{g}$.

\section{Generic Matrices}

In this section we define generic matrices and construct a relatively free algebra in the class determined by $M_n(K)$ with an elementary grading induced by $\textbf{g}=(g_1,\dots,g_n)\in G^n$. We also prove some results that will be used in the next sections. 

For each $h \in G$ let $Y_h=\{y_{h,i}^{k},  1\leq k \leq n, i = 1,2,\dots \}$ be a countable set of commuting variables. Let $Y=\cup_{h \in G}Y_h$ and denote by $\Omega=K[Y]$ the commutative polynomial algebra generated by $Y$. The set $\{g_1,\dots,g_n\}$ is denoted by $G_n$. 

The algebra $M_n(\Omega)$ has an elementary grading induced by the $n$-tuple $\textbf{g}=(g_1,\dots,g_n)$, given $h\in G$ we are interested in determining the elementary matrices $E_{ij}$ of degree $h$. If we fix $i \in I_n$ there exists an elementary matrix $E_{ij}$ of degree $h$ if and only if $g_ih \in G_n$ and in this case $j\in I_n$ is determined by the equality $g_j=g_ih$. 

For each $h \in G$ we denote by $L_{h}$ the set of all indexes $k \in I_n$ such that $g_kh \in G_n$ and by $s^k_h\in I_n$ is the index determined by $g_kh=g_{s^k_h}$. Then it is easy to see that $(M_n(\Omega))_h=0$ iff $L_h=\emptyset$, moreover if $L_h \neq \emptyset$ then $\{E_{ks_h^k}|k\in L_h\}$ is the set of all elementary matrices of degree $h$. We consider in $M_n(\Omega)$ the homogeneous matrices
\begin{equation}\label{e3}
A_i^{(h)}=\sum_{k \in L_{h}} y_{h,i}^{k}E_{k,s^k_h}.
\end{equation}
These matrices will be called \textsl{generic matrices}, the subalgebra $F$ of $M_n(\Omega)$ generated by the generic matrices is a graded subalgebra, in this grading the generic matrix $A_i^{h}$ defined above is homogeneous of degree $h$. 

\begin{lemma}\label{l4}
The relatively free algebra $K\langle X \rangle / T_G(M_n(K))$ is isomorphic to the algebra $F$.
\end{lemma}

\textit{Proof.} The proof is analogous to that of \cite[Lemma 1]{Azevedo2}. \hfill $\Box$

\begin{remark}
As a direct consequence we have $T_G(F)=T_G(M_n(K))$, so from now on we will focus on the graded identities for the algebra $F$.
\end{remark}

In order to proceed we need convenient notation to compute a product of generic matrices. 

\begin{definition}\label{definition}
Let $\textbf{h}=(h_1,\dots,h_q)\in G^{q}$ , the set \[L_{\textbf{h}}=\{k\in I_n \mid g_kh_1\dots h_i\in G_n, \mbox{\textrm{ for every }} i \in I_q\}\] is the \textbf{set associated with h}. For each $k \in L_{\textbf{h}}$ define the $(q+1)$-tuple $s_k=(s_1^k,\dots,s_q^k,s_{q+1}^k)$, inductively by setting:

\begin{itemize}
\item[(1)] $s_1^{k}=k$,
\item[(2)] for $i\in I_{q}$ the index $s_{k+1} \in I_n$ is determined by $g_{s_{i+1}^k}=g_{s_{i}^k}h_{i}$.
\end{itemize}
\end{definition}

\begin{remark}
In the above definition $k \in L_{\textbf{h}}$ if and only if there exist elementary matrices $E_{i_1j_1},\dots E_{i_qj_q}$ such that $i_1=k$, $E_{i_aj_a}$ has degree $h_a$ and $E_{i_1j_1}\cdot \dots E_{i_qj_q}\neq 0$. Moreover $i_a=s^{k}_a$ and $j_a=s^{k}_{a+1}$.
\end{remark}

\begin{lemma}\label{l1}
If $L$ is the set of indexes associated with the $q$-tuple $(h_1,\dots,h_q)$ in $G^{q}$ and $s_k=(s_1^{k},\dots,s_q^{k},s_{q+1}^k)$ denotes the corresponding sequence determined by $k\in L$ then 
\[A_{i_1}^{h_1}\dots A_{i_q}^{h_q}=\sum_{k \in L}w_kE_{s_1^{k},s_{q+1}^k}\]
where $w_k=y_{h_1,i_1}^{s_{1}^{k}}y_{h_2,i_2}^{s_{2}^{k}}\dots y_{h_q,i_q}^{s_{q}^{k}}$.
\end{lemma}

\textit{Proof.} From the previous remark we conclude that \[E_{k_1s_{h_1}^{k_1}}E_{k_2s_{h_2}^{k_2}}\dots E_{k_qs_{h_q}^{k_q}}\neq 0,\] iff $k_1 \in L$ and for every $i \in I_q$ we have $k_i=s_i^k$, $s_{h_1}^{k_i}=s_{i+1}^{k}$. Since from (\ref{e3}) we have 
\[A_{i_1}^{h_1}\dots A_{i_q}^{h_q}=\sum (y_{h_1,i_1}^{k_1}\dots y_{h_q,i_q}^{k_q})E_{k_1,s^{k_1}_{h_1}}\dots E_{k_q,s^{k_q}_{h_q}},\]the result follows.
\hfill $\Box$

\begin{remark}
In order to simplify the notation we will adopt the following convention. If $f(x_{i_1},\dots,x_{i_k})\in K\langle X \rangle$ is a polynomial in the variables $x_{i_1},\dots,x_{i_k}$ then $A_{i_l}$ will denote the generic matrix $A_{i_l}^{h(x_{i_l})}$, and $f(A_{i_1}\dots A_{i_k})$ denotes the result of substituting each variable for the corresponding generic matrix.
\end{remark}
The following two consequences of the above lemma will be useful in the next section. We recall that given a monomial $m=x_{i_1}\dots x_{i_q}$ of length $q$, the sequence $(h_1,\dots,h_q)$ where $h_k$ is the degree of the variable $x_{i_k}$, is denoted by $h(m)$ and $\alpha(m)$ denotes the degree of the monomial in the $G$-grading of $K\langle X \rangle$ defined in Section \ref{preliminaries}.

\begin{corollary}\label{c1}
Let $m_1,m_2$ be monomials such that $h(m_1)=h(m_2)$, then $m_1\in T_G(F)$ if and only if $m_2 \in T_G(F)$.
\end{corollary}

\textit{Proof.} It follows directly from Lemma \ref{l1} above since $m_i \in T_G(F)$ if and only if the set associated with $h(m_1)$ is empty.
\hfill $\Box$

\begin{corollary}\label{c2}
If $x_{i_1}x_{i_2}\dots x_{i_r}$ and $x_{j_1}x_{j_2}\dots x_{j_s}$ are two monomials such that the matrices $A_{i_1}\dots A_{i_r}$ and $A_{j_1}\dots A_{j_s}$ have in the same position the same non-zero entry then $r=s$ and there exists $\sigma \in S_r$ such that 
\[x_{j_1}x_{j_2}\dots x_{j_s}=x_{i_{\sigma(1)}}x_{i_{\sigma(2)}}\dots x_{i_{\sigma(r)}},\] 
and moreover for every $l\in I_r$ \[\alpha(x_{i_{\sigma(1)}}\dots x_{i_{\sigma(l-1)}})=\alpha(x_{i_1}\dots x_{i_{\sigma(l)-1}}).\]
\end{corollary}

\textit{Proof.} Let us assume the matrices $A_{i_1}\dots A_{i_r}$ and $A_{j_1}\dots A_{j_s}$ have the same non-zero entry in the position $(k,l)$, let $\textbf{h}(m)=(h_1^{m},\dots, h_r^m)$ and $\textbf{h}(n)=(h_1^{n},\dots, h_s^n)$. It follows from Lemma \ref{l1} that $k \in L_{\textbf{h}(m)}\cap L_{\textbf{h}(n)}$, moreover the two monomials in $\Omega$ \[y_{h_1^{m},i_1}^{s_{1}^{k,m}}\cdot y_{h_2^{m},i_2}^{s_{2}^{k,m}}\cdot \dots \cdot  y_{h_r^{m},i_r}^{s_{r}^{k,m}}\] and \[y_{h_1^{n},j_1}^{s_{1}^{k,n}}\cdot y_{h_2^{n},j_2}^{s_{2}^{k,n}}\cdot \dots \cdot y_{h_s^{n},j_s}^{s_{s}^{k,n}}\] are equal, where $s^{k,m}=(s^{k,m}_1,s^{k,m}_2,\dots,s^{k,m}_r)$ (resp. $s^{k,n}$) denotes de sequence corresponding to $k\in L_{\textbf{h}(m)}$ (resp. $k \in L_{\textbf{h}(n)}$). 

From the equality of the monomials we conclude that $r=s$ and there exists $\sigma \in S_r$ such that $j_l=i_{\sigma(l)}$ for all $l \in I_r$, then \[x_{j_1}x_{j_2}\dots x_{j_r}=x_{i_{\sigma(1)}}x_{i_{\sigma(2)}}\dots x_{i_{\sigma(r)}}.\] Moreover $s_l^{k,n}=s_{\sigma(l)}^{k,m}$, and it follows from Definition \ref{definition} that \[g_kh_{i_{\sigma(1)}}\dots h_{i_{\sigma(l-1)}}=g_kh_{i_1}\dots h_{i_{\sigma(l)-1}},\]  therefore $\alpha(x_{i_{\sigma(1)}}\dots x_{i_{\sigma(l-1)}})=\alpha(x_{i_1}\dots x_{i_{\sigma(l)-1}})$. 
\hfill $\Box$

\section{Graded Identities and Preliminary Results}

We consider the following polynomials:

\begin{equation}\label{1}
x_1x_2-x_2x_1,\mbox{ }h(x_1)=h(x_2)=\epsilon,
\end{equation}

\begin{equation}\label{2}
x_{1}x_{3}x_{2}-x_{2}x_{3}x_{1},\mbox{ }\epsilon \neq h(x_1)=h(x_2)=h(x_3)^{-1},
\end{equation} 

\begin{equation}\label{3}
x_1=0\mbox{ }, (M_n(K))_{h(x_1)}=0.
\end{equation}

\begin{proposition}\label{p1}
Let $\textbf{g}=(g_1,\dots, g_n)\in G^{n}$ be an arbitrary $n$-tuple and let $M_n(K)=\bigoplus_{g \in G}(M_n(K))_g$ be the elementary grading induced by $\textbf{g}$. The following statements are equivalent:
\begin{enumerate}
\item[(i)] If $i\neq j$ then $g_i \neq g_j$, i. e., the elements in $\textbf{g}$ are pairwise different;
\item[(ii)] The subspace $(M_n(K))_{\epsilon}$ coincides with the subspace of the diagonal matrices;
\item[(iii)] The polynomial $x_1x_2-x_2x_1$, where $h(x_1)=h(x_2)=\epsilon$, is a graded identity for $M_n(K)$.
\end{enumerate}
\end{proposition}

\textit{Proof.}
Clearly $(i)$ and $(ii)$ are equivalent and $(ii)$ implies $(iii)$. To conclude we will prove that $(iii)$ implies $(i)$. Let $E_{ij}$ be an elementary matrix of degree $\epsilon$, it follows from $(iii)$  that $E_{ij}=E_{ii}E_{ij}=E_{ij}E_{ii}$, hence $i=j$.
\hfill $\Box$

\begin{remark}\label{r1}
If there are $m$ equal elements in the $n$-tuple it is easy to see that there is a subalgebra of $(M_n(K))_{\epsilon}$ isomorphic to $M_{m}(K)$ and any ordinary identity of this algebra would produce a graded identity of $M_n(K)$ in the variables $X_{\epsilon}$.
\end{remark}

We recall that, with the exception of Proposition \ref{p1} and Remark \ref{r1} above, we fixed an $n$-tuple $\textbf{g}=(g_1,\dots, g_n)\in G^{n}$ of pairwise different elements and $M_n(K)=\bigoplus_{g \in G}(M_n(K))_g$ is the elementary grading induced by $\textbf{g}$.

\begin{lemma}\label{l3}
The $G$-graded algebra $M_n(K)$ satisfies the $G$-graded polynomial identities (\ref{1}), (\ref{2}) and (\ref{3}).
\end{lemma}

\textit{Proof.} \cite[Lemma 4.1]{drenskyebahturin}. \hfill $\Box$

\begin{definition}
We denote by $J$ the $T_G$-ideal generated by the identities (\ref{1}), (\ref{2}) and (\ref{3}).
\end{definition}

\begin{lemma}\label{l2}
Let $\overline{m}(x_1,\dots,x_q)$ and $\overline{n}(x_1,\dots,x_q)$ be two monomials that start with the same variable and let $m(x_1,\dots, x_q)$, $n(x_1,\dots, x_q)$ be the monomials obtained from $\overline{m}$ and $\overline{n}$ respectively by deleting the first variable. If there exist matrices $A_1,\dots, A_q$, such that $\overline{m}(A_1,\dots,A_q)$ and $\overline{n}(A_1,\dots,A_q)$ have in the same position the same non-zero entry then $m(A_1,\dots,A_q)$ and $n(A_1,\dots,A_q)$ also have in the same position the same non-zero entry.
\end{lemma}
\textit{Proof.} It follows directly from Lemma \ref{l1}.
\hfill $\Box$

In the next lemma we follow the idea of Azevedo, see \cite[Lemma 6]{Azevedo1} and \cite[Lemma 5]{Azevedo2}.

\begin{lemma}\label{l5}
Let $m(x_1,\dots,x_p)$ and $n(x_1,\dots,x_p)$ be two monomials such that the matrices $n(A_1,\dots,A_p)$ and $m(A_1,\dots,A_p)$ have in the same position the same non-zero entry then \[m(x_1,x_2\dots,x_p)\equiv n(x_1,x_2\dots,x_p) \mbox{ modulo }J.\]
\end{lemma}

\textit{Proof.} 
Let $m(x_1,\dots,x_p)=x_{i_1}\dots x_{i_q}$, for integers $1\leq k <l \leq q+1$ we denote by $m^{[k,l]}=x_{i_{k}}\dots x_{i_{l-1}}$ (resp. $n^{[k,l]}$) the monomial obtained from $m$ (resp. $n$) by deleting the first $k-1$ variables and the last $q-l+1$ variables. 

It follows from Corollary \ref{c2} that there exists $\sigma \in S_q$ such that
\begin{equation}\label{e2}
 n(x_1,\dots,x_p)=x_{i_{\sigma(1)}}\dots x_{i_{\sigma(q)}},
\end{equation}
  and 
\begin{equation}\label{e1}
 \alpha(n^{[1,l]})=\alpha(m^{[1,\sigma(l)]}),
\end{equation}
for all $l \in I_{q}$. We will prove the result by induction on $q$, if $q=1$ the lemma is clearly true. 

Let $a=\sigma^{-1}(1)$, it follows from (\ref{e1}) that \[\alpha(n^{[1,a]})=\epsilon.\] Assume there exists $r \in I_{q-1}$ such that \[\sigma^{-1}(r)<a<\sigma^{-1}(r+1).\] The result follows from the previous lemma and the induction hypothesis if we prove that in this case $n(x_1,\dots,x_p)$ is congruent modulo $J$ to a monomial that starts with $x_{i_1}$.
It follows from (\ref{e1}) that \[\alpha(n^{[\sigma^{-1}(r),\sigma^{-1}(r+1)]})=\alpha(m^{[r,r+1]})=\alpha(x_{i_r}).\] Using (\ref{e2}) we obtain $n^{[\sigma^{-1}(r),\sigma^{-1}(r+1)]}=x_{i_r}n^{[\sigma^{-1}(r)+1,\sigma^{-1}(r+1)]}$, therefore \[\alpha(n^{[\sigma^{-1}(r)+1,\sigma^{-1}(r+1)]})=\epsilon.\] Finally if $\sigma^{-1}(r)+1=a$ then using the identity (\ref{1}) we have \[n=n^{[1,a]}n^{[a,\sigma^{-1}(r+1)]}\equiv_J n^{[a,\sigma^{-1}(r+1)]}n^{[1,a]}=x_{i_1}n^{[a+1,\sigma^{-1}(r+1)]}n^{[1,a]}.\] The other possibility is $\sigma^{-1}(r)+1<a$, in this case we have \[n=n^{[1,\sigma^{-1}(r)+1]}n^{[\sigma^{-1}(r)+1,a]}n^{[a,\sigma^{-1}(r+1)]},\] since \[\alpha(n^{[1,\sigma^{-1}(r)+1]}n^{[\sigma^{-1}(r)+1,a]})=\alpha(n^{[\sigma^{-1}(r)+1,a]}n^{[a,\sigma^{-1}(r+1)]})=\epsilon\] it follows from (\ref{2}) that $n\equiv_Jn^{[a,\sigma^{-1}(r+1)]}n^{[\sigma^{-1}(r)+1,a]}n^{[1,\sigma^{-1}(r)+1]}$, and the first variable in this last monomial is $x_{i_1}$. Hence this lemma is proved under the assumption that exists $r \in I_{q-1}$ such that $\sigma^{-1}(r)<a<\sigma^{-1}(r+1)$.

If there exists no such $r$ it follows that $\sigma(I_{a-1})=\{b,b+1,\dots,q\}$ for some $b \in I_q$, and in this case the monomials $m^{[b,q+1]}$ and $n^{[1,a]}$ have the same multidegree and in particular $\alpha (m^{[b,q+1]})=\alpha(n^{[1,a])}=\epsilon$. Moreover $\sigma(1)\geq b$ and it follows from (\ref{e1}) that $\alpha(m^{[1,\sigma(1)]})=\alpha(n^{[1,1]})=\epsilon$, hence \[\alpha(m^{[1,b]}m^{[b,\sigma(1)]})=\alpha(m^{[b,\sigma(1)]}m^{[\sigma(1),q+1]})=\epsilon\] and since $m=m^{[1,b]}m^{[b,\sigma(1)]}m^{[\sigma(1),q+1]}$, it follows from $(2)$ that \[m\equiv_J m^{[\sigma(1),q+1]}m^{[b,\sigma(1)]}m^{[1,b]},\] and this last monomial starts with the same variable as $n$.
\hfill $\Box$

\section{A basis for the graded identities of $M_n(K)$}

In \cite{Azevedo1, Azevedo2} it is shown that there are no nontrivial identities $x_{i_1}\dots x_{i_k}$ for $M_n(K)$ with its natural elementary $\mathbb{Z}_n$ and $\mathbb{Z}$ gradings respectively, however in \cite[Example 4.7, Theorem 4.8, Theorem 4.9]{drenskyebahturin} concrete gradings are considered where $M_n(K)$ satisfies nontrivial identities $x_{i_1}\dots x_{i_k}$. These results also hold for infinite fields of positive characteristic. 

In this section in Theorem \ref{t1} we prove that the main result of \cite{drenskyebahturin} holds for arbitrary infinite fields. T		his theorem states that for the elementary gradings considered here a basis for the graded identities of $M_n(K)$ consists of (\ref{1})-(\ref{3}), and a finite number of monomials $x_{i_1}\dots x_{i_k}$ of length $k$ bounded by a function of $n$. The following lemma corresponds to \cite[Proposition 4.2]{drenskyebahturin}.

\begin{lemma}\label{previous}
Let $G_0 = \{g \in G | (M_n(K))_g \neq 0\}$ be the support of a $G$-grading
of $M_n(K)$ and let $I$ be the set of all finite sequences $(h_1,\dots, h_k)$ of elements
of $G_0$ such that the monomial $m$ with $h(m)=(h_1,\dots, h_k)$ is a G-graded polynomial identity of $M_n(K)$. Then there exists a positive integer $n_0$ such that $m$ is a consequence
of the $G$-graded polynomial identities of $M_n(K)$ as in (\ref{1}), (\ref{2}) and (\ref{3}) together with the monomials $m$ of length $k$, where $h(m)\in I$ and $k<n_0$.
\end{lemma}

\textit{Proof} Let $s=|G_0|$, $n_0=4s^{2s+2}$ and let $U$ be the $T_G$-ideal generated by (\ref{1})-(\ref{3}) and the monomials $m$ of length less then $n_0$ such that $h(m)\in I$. If $m \in T_G(F)$ is a multilinear monomial then it follows from \cite[Proposition 4.2]{drenskyebahturin} that $m \in U$. If $m$ is not multilinear let $\overline{m}$ be a multilinear monomial such that $h(m)=h(\overline{m})$. It follows from Corollary \ref{c1} that $\overline{m}\in T_G(F)$, then $\overline{m}$ is in $U$ and since $m$ is obtained from $\overline{m}$ by identifying some of the variables we conclude that $m \in U$.
\hfill $\Box$

\begin{theorem}\label{t1}
Let $K$ be an infinite field. Let $G$ be any group and let $g = (g_1,\dots, g_n)\in G^n$ induce an elementary $G$-grading of $M_n(K)$ where the elements $g_1,\dots,g_n$ are pairwise different.
Then a basis of the graded polynomial identities of $M_n(K)$ consists of (\ref{1})-(\ref{3})
and a finite number of identities of the form $x_{i_1}\dots x_{i_k}$, $k\geq 2$, where the
length $k$ is bounded by a function of $n$.
\end{theorem}

\textit{Proof.} Let $U$ be the $T_G$-ideal defined in the previous lemma, since $|G_0|=s\leq n^2$ it is easy to see that we may choose $k\leq 4n^{4(n^2+1)}$. It follows from Lemma \ref{l3} that $U\subset T_G(F)$. We wish to prove that $T_G(F)\subset U$, suppose on contrary that the inclusion does not hold. Since the field $K$ is infinite, there exists a multihomogeneous identity of $F$ that does not belong to $U$. Let \[f(x_1,\dots,x_n)=\sum_{i=1}^{k_0}a_i m_i(x_1,\dots,x_n),\] where $a_i \neq 0$, $1\leq i \leq k_0$, be a multihomogeneous element of $T_G(F)-U$ with the minimal number of non-zero summands. If $k_0=1$ then the monomial $m_1$ is an identity of $F$ that is not an element of the $T_G$-ideal $U$, but this contradicts Lemma \ref{previous}, hence $k_0>1$. Moreover, since $k_0$ is minimal we conclude that $m_1$ is not an identity. Clearly, \[-a_1m_1(A_1,\dots,A_n)=\sum_{i=2}^{k_0}a_im_i(A_1,\dots,A_n).\] The non-zero entries in each matrix $m_i(A_1,\dots,A_n)$ are monomials in $\Omega$ and $m_1(A_1,\dots,A_n)\neq 0$, hence there exists $p\in \{2,\dots, k_0\}$ such that the matrices $m_1(A_1,\dots,A_n)$ and $m_p(A_1,\dots,A_n)$ have in the same position the same nonzero entry. Using Lemma \ref{l5} we conclude that $m_p - m_1 \in U$. In this case \[f(x_1,\dots,x_n)+a_p(m_1(x_1,\dots,x_n)-m_p(x_1,\dots,x_n))\in T_G(F)-U,\] but the last polynomials is a sum of $k_0-1$ monomials and this is a contradiction since $k_0$ is minimal.
\hfill $\Box$

The following corollary generalizes the main result of \cite{Azevedo1} for any elementary grading of $M_n(K)$ in which the neutral component consists of the diagonal matrices and $|G|=n$. 

\begin{corollary}
Let $K$ be an infinite field. Let $G$ be a finite group of order $n$ and let $g = (g_1,\dots, g_n)\in G^n$ induce an elementary $G$-grading of $M_n(K)$ where the elements $g_1,\dots,g_n$ are pairwise different.
Then a basis of the graded polynomial identities of $M_n(K)$ consists of (\ref{1}) and (\ref{2}).
\end{corollary}
\textit{Proof.} Since $|G|=n$ the set $G_n$ in Definition \ref{definition} is the group $G$, therefore for every $q\in \mathbb{N}$ and every sequence $\textbf{h} \in G^{q}$ the set $L_{\textbf{h}}$ associated with this sequence is $I_n=\{1,2,\dots, n\}$. Hence Lemma \ref{l1} implies no monomial $x_{i_1}\dots x_{i_k}\in K\langle X \rangle$ of length $k\geq 1$ is a graded identity for $M_n(K)$ and the corollary now follows from the previous theorem.
\hfill $\Box$

\begin{flushleft}
\textbf{Acknowledgements}
\end{flushleft}
We thank the referee for valuable remarks. Thanks are due to P. Koshlukov for his useful suggestions.


\begin{thebibliography}{99}

\bibitem {Azevedo1}\textrm{S. S. Azevedo},
\textit{Graded identities for the matrix algebra of order n over an infinite field}, Comm. Algebra \textbf{30}, 5849--5860, 12 (2002).

\bibitem {Azevedo2}\textrm{S. S. Azevedo},
\textit{A basis for $\mathbb{Z}$-graded identities of matrices over infinite fields}, Serdica Math. Journal \textbf{29 (2)}, 149--158 (2003).

\bibitem {AzevedoeKoshlukov} \textrm{S. S. Azevedo, P. Koshlukov},
\textit{Graded identities for T-prime algebras over field of positive characteristic}, Israel J. Math. \textbf{128}, 157--176 (2002).

\bibitem{drenskyebahturin} \textrm{Y. Bahturin, V. Drensky},
\textit{Graded polynomial identities of matrices}, Linear Algebra and its Applications, 15--34 (2002).

\bibitem {Bahturin} \textrm{Yu. A. Bakhturin, M. V. Zaicev}, \textit{Group Gradings on Matrix Algebras}, Candad. Math. Bull. \textbf{45 (4)}, 499--508 (2002).

\bibitem {Berele} \textrm{ A. Berele}, \textit{Magnum PI}, Israel J. Math. \textbf{51}, no. 1-2, 13--19, (1985).

\bibitem {BelovRowen} \textrm{A. Ya. Belov, L. H. Rowen},
\textit{Computational Aspects of polynomial identities}, Research Notes in Mathematics \textbf{9}, A.K. Peters, Ltd., Wellesley, MA, 2005.

\bibitem {Drensky} \textrm{ V. Drensky},
\textit{A minimal basis of identities for a second-order matrix algebra over a  field of characteristic 0}, Algebra and Logic \textbf{20(3)}, 188--194(1981).

\bibitem {drensky}\textrm{V. Drensky},
\textit{Free algebras and PI algebras}, Graduate Course in Algebra, Springer-Verlag PTE.LTD, (1999).

\bibitem{drenskyformanek} \textrm{V. Drensky, E. Formanek},
\textit{Polynomial Identity Rings}, CRM Advanced Courses in Mathematics, Birkhäuser, Basel (2004).

\bibitem{formanek} \textrm{E. Formanek},
\textit{The ring of generic matrices}, J. Algebra, \textbf{258}, no. 1, 310--320 (2002).

\bibitem {Genov} \textrm{G. K. Genov},
\textit{A basis for identities of third order matrix algebra over a finite field}, Algebra and Logic, \textbf{20}, 241--257(1981).

\bibitem{GenoveSiderov}\textrm{G. K. Genov, P.N. Siderov},
\textit{A basis for identities of the algebra of fourth-order matrices over a finite field I}, II Serdica 8, 313--323, 351--366(1982).

\bibitem {GiambrunoZaicev} \textrm{A. Giambruno, M. Zaicev},
\textit{Polynomial identities and asymptotic methods}, Math. Surv. and Monographs \textbf{122}, AMS.

\bibitem {IoneMontes} \textrm{S. D$\check{a}$sc$\check{a}$lescu, B. Ion, C. N$\check{a}$st$\check{a}$sescu, J. Rios Montes},
\textit{Group gradings on full matrix rings}, J. Algebra \textbf{220}, 709--728, 1999.

\bibitem {Jacobson} \textrm{N. Jacobson}, \textit{PI-Algebras, an introduction}, Springer Lecture Notes in Math. \textbf{441}, Springer, (1975).

\bibitem {Kaplansky} \textrm{I. Kaplansky}, \textit{Rings with a polynomial identity}, Bull. Amer. Math. Soc. \textbf{54} \textbf{220}, 496--500, 1948.

\bibitem {Kemer} \textrm{A. Kemer},
\textit{Finite basis property of identities of associative algebras}, Algebra and Logic \textbf{26}, 362--397(1987).

\bibitem {Koshlukov} \textrm{P. Koshlukov},
\textit{Basis of the identities of the matrix algebra of order two over a field of caracteristic $p\neq 2$}, J. Algebra \textbf{241}, 410--434 (2001).

\bibitem {MaltseveKuzmin} \textrm{Yu. N. Maltsev, E. N. Kuzmin},
\textit{A basis for the identities of the algebra of second-order over a finite field}, Algebra and Logic \textbf{17} No. 1, 18--21(1978).

\bibitem {Razmyslov} \textrm{Yu. P. Razmyslov},
\textit{Finite basing of the identities of a matrix algebra of second order over a field of characteristic zero},  Algebra and Logic \textbf{12}, 47--63(1973).

\bibitem {Razmyslov2} \textrm{Yu. P. Razmyslov},
\textit{Trace Identities of full matrix algebras over a field of characteristic zero},  Izv. Akad. Nauk SSSR, Ser. Mat. \textbf{38}, 723--756 (1974). 

\bibitem {Rowen} \textrm{L. H. Rowen},
\textit{Polynomial Identities in Ring Theory}, Acad. Press Pure and Applied Math., vol. 84, New York (1980).

\bibitem {Vasilovsky1} \textrm{S. Yu. Vasilovsky},
\textit{$\mathbb{Z}$-graded polynomial identities of the full matrix algebra}, Commun. Algebra \textbf{26 (2)},  601--612 (1998).

\bibitem {Vasilovsky2} \textrm{S. Yu. Vasilovsky},
\textit{$\mathbb{Z}_n$-graded polynomial identities of the full matrix algebra of order $n$}, Proc. Amer. Math. Soc.  \textbf{127 (12)},  3517--3524 (1999).

\bibitem {Vincenzo} \textrm{O. M. Di Vincenzo},
\textit{On the graded identities of $M_{1,1}(E)$}, Isr. J. Math. \textbf{80}, 323--335 (1992).

\bibitem {GiambrunoKoshlukov} \textrm{O. M. Di Vincenzo},
\textit{On the graded identities of $M_{1,1}(E)$}, Isr. J. Math. \textbf{80}, 323--335 (1992).

\end{thebibliography}
\end{document}